\documentclass[12pt]{article}
\usepackage{amssymb,latexsym,amsmath}
\usepackage{amsfonts}
\usepackage[cp1251]{inputenc}
\usepackage[epsf]{}
\usepackage[english,russian]{babel}
\usepackage[dvips]{graphicx}
\newtheorem{thm}{Theorem}

\newtheorem{lm}{Lemma}
\newtheorem{cor}{Corollary}
\newtheorem{step}{Step}

\begin{document}

\renewcommand{\refname}{References}
\renewcommand{\abstractname}{Absrtact}
\renewcommand{\figurename}{Fig.}

\title{On spectral decomposition  of Smale-Vietoris axiom A diffeomorphisms
\footnote{2000{\it Mathematics Subject Classification}. Primary 37D20; Secondary 37G30}
}

\author{N.~Isaenkova$^{1}$\and E.~Zhuzhoma$^{2}$}
\date{}
\maketitle

{\small $^{1}$ Nizhny Novgorod Academy MVD of Russia, RUSSIA, nisaenkova@mail.ru}\\
{\small $^{2}$ National Research University Higher School of Economics, RUSSIA, zhuzhoma@mail.ru}

\begin{abstract}
We introduce Smale-Vietoris diffeomorphisms that include the classical DE-mappings with Smale solenoids. We describe the correspondence between basic sets of axiom A Smale-Vietoris diffeomorphisms and basic sets of nonsingular axiom A endomorphisms. For Smale-Vietoris diffeomorphisms of 3-manifolds, we prove the uniqueness of nontrivial solenoidal basic set. We construct a bifurcation between different types of solenoidal basic sets which can be considered as a destruction (or birth) of Smale solenoid.
\end{abstract}

\section*{Introduction}

\begin{figure}[h]                                                                       
\begin{center}\includegraphics[width=14cm,height=3cm]{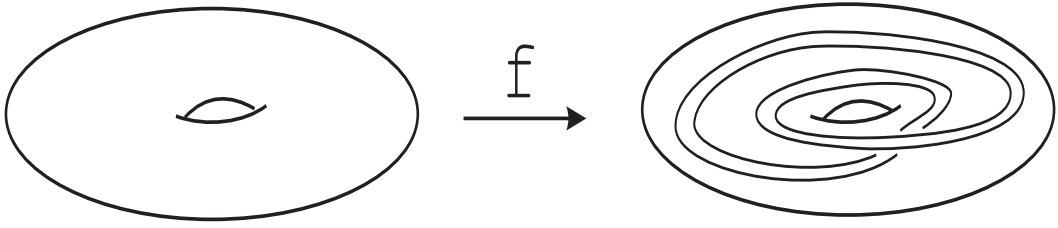}\end{center}
\caption{DE-map by S. Smale} \label{SmaleSol}
\end{figure}

Stephen Smale, in his celebrated paper \cite{Smale67}, introduced so-called DE-maps which arise from expanding maps (the abbreviation DE is formed by first letters of \textit{Derived} from \textit{Expanding} map). Let $T$ be a closed manifold of dimension at least 1, and $N$ an $n$-disk of dimension $n\geq 2$. Omitting details, one can say that a DE map is the skew map
\begin{equation}\label{eq:map}
 f: T\times N\to T\times N,\quad (x;y)\mapsto\left(g_1(x);\,\, g_2(x,y)\right),
 \end{equation}
where $g_1: T\to T$ is an expanding map of degree $d\geq 2$, and
 $$ g_2|_{\{x\}\times N}: \{x\}\times N\to \{g_1(x)\}\times N $$
an uniformly attracting map of $n$-disk $\{x\}\times N$ into $n$-disk $\{g_1(x)\}\times N$ for every $x\in T$. In addition, $f$ must be a diffeomorphism onto its image $T\times N\to f(T\times N)$. In the particular case, when $T=S^1$ is a circle and $N=D^2$ is a 2-disk with the uniformly attracting $g_2$, one gets a classical Smale solenoid $\cap_{l\ge 0}f(\mathbb{T}\times D^2)=\mathfrak{S}$, see Fig.~\ref{SmaleSol}, that is a topological solenoid.

Recall that a topological solenoid was introduced by Vietoris \cite{Vietoris27} in 1927 (independently, a solenoid was introduced by van Danzig \cite{vanDanzig1930} in 1930, see review in \cite{Takens2005}). Smale \cite{Smale67} proved that $\mathfrak{S}(f)$ is a hyperbolic expanding attractor. This construction was generalized by Williams \cite{Williams67,Williams74} who defined $g_1$ to be expansion mappings of branch manifolds (this allows to Williams to classify interior dynamics of expanding attractors) and by Block \cite{BlockL75}
who defined $g_1$ to be an axiom A endomorphism. The last paper concerns to the $\Omega$-stability and the proving of decomposition of non-wandering set into so-called basic sets (Spectral Decomposition Theorem for A-endomorphisms). Ideologically, our paper is a continuation of \cite{BlockL75}, where it was proved the following result (Theorem A).
Let $f: M^n\to M^n$ be a Smale-Vietoris diffeomorphism of closed $n$-manifold $M^n$ and $\mathfrak{B}\subset M^n$ the support of Smale skew-mapping $f|_{\mathfrak{B}}$ (see the notations below). Then $f|_{\mathfrak{B}}$ satisfies axiom A on $\mathfrak{B}$ if and only if $g$ does on $T$.

Let us mention that in the frame of Smale-Williams construction the interesting examples of expanding attractors was obtained in \cite{Bothe92,FarrellJones80,IsaenkovaZhuzhoma2009,Jones83,RobinsonWilliams76}. Bothe \cite{Bothe83} classified the purely Smale solenoids on 3-manifolds. He was first who also proved that a DE map $S^1\times D^2\to S^1\times D^2$ can be extended to a diffeomorphism of some closed 3-manifold $M^3\supset S^1\times D^2$ (see also \cite{JiagNiWang2004,JimingMaBinYu2007,JimingMaBinYu2011}). Ya. Zeldovich and others (see \cite{ChildressGilbert-book-1995}) conjectured that Smale type mappings could be responsible for so-called fast dynamos. Therefore, it is natural to consider various generalizations of classical Smale mapping.

In a spirit of Smale construction of DE-maps, we here introduce diffeomorphisms called Smale-Vietoris that are derived from nonsingular endomorphisms. A non-wandering set of Smale-Vietoris diffeomorphisms belongs to an attractive invariant set of solenoidal type. In the classical case, the invariant set coincides with the non-wandering set consisting of a unique basic set. In general, the non-wandering set does not coincide with the invariant set, and divides into basic sets provided the nonsingular endomorphism is an A-endomorphism.

Let $N$ be $(n-k)$-dimensional compact Riemannian manifold with a non-empty boundary where $n-k\geq 1$. For a subset $N_1\subset N$, we define the diameter $diam~N_1=\max_{a,b\in N_1}\{\rho_N(a,b)\}$ of $N_1$ where $\rho_N$ is the metric on $N$. Denote by $\mathbb{T}^k=\underbrace{S^1\times\cdots\times S^1}_{k}$ the $k$-dimensional torus, $k\in\mathbb{N}$.  A surjective mapping $g: \mathbb{T}^k\to \mathbb{T}^k$ is called a $d$-\textit{cover} if $g$ is a preserving orientation local homeomorphism of degree $d$. A good example is the preserving orientation linear expanding mapping $E_d: \mathbb{T}^k\to \mathbb{T}^k$ defined by an integer $k\times k$ matrix with the determinant equals $d$. Certainly, $E_d$ is a $d$-cover.

A skew-mapping
\begin{equation}\label{eq:degree-d}
  F: \mathbb{T}^k\times N\to \mathbb{T}^k\times N,\quad (t,z)\longmapsto\left((g(t);\,\,\omega(t,z)\right)
\end{equation}
is called a \textit{Smale skew-mapping} is the following conditions hold:
 \begin{itemize}
    \item $F: \mathbb{T}^k\times N\to F\left(\mathbb{T}^k\times N\right)$ is a diffeomorphism on its image;
    \item $g: \mathbb{T}^k\to \mathbb{T}^k$ is a $d$-cover, $d\geq 2$;
    \item given any $t\in \mathbb{T}^k$, the restriction $w|_{\{t\}\times N}: \{t\}\times N\to \mathbb{T}^k\times N$
          is the uniformly attracting embedding
    \begin{equation}\label{eq:pavnom}
    \{t\}\times N\to int\left(\{g(t)\}\times N\right)
    \end{equation}
          i.e., there are $0<\lambda <1$, $C>0$ such that
    \begin{equation}\label{eq:contraction}
    diam~(F^n(\{t\}\times N))\le C\lambda^n diam~(\{t\}\times N),\,\,\forall\, n\in\mathbb{N}.
    \end{equation}
 \end{itemize}
When $g=E_d$, Smale skew-mapping is a $DE$ mapping (\ref{eq:map}) introduced by Smale \cite{Smale67}.

A diffeomorphism $f: M^n\to M^n$ is called a \textit{Smale-Vietoris} diffeomorphism if there is the $n$-submanifold
$\mathbb{T}^k\times N \subset M^n$ such that the restriction $f|_{\mathbb{T}^k\times N}\stackrel{\rm def}{=}F$ is a Smale skew-mapping. The submanifold
$\mathbb{T}^k\times N \subset M^n$ is called a \textit{support of Smale skew-mapping}.

Put by definition,
 $$ \cap_{l\ge 0}F^l(\mathbb{T}^k\times N)\stackrel{\rm def}{=}\mathfrak{S}(f). $$
One can easy to see that the set $\mathfrak{S}(f)=\mathfrak{S}$ is attractive, invariant and closed, so that the restriction
 $$ f|_{\mathfrak{S}}: \mathfrak{S}\to\mathfrak{S} $$
is a homeomorphism.

The following theorem shows that there is an intimate correspondens between basic sets of $f|_{\mathfrak{B}}$ and basic sets of the A-endomorphism $g$.

\begin{thm}\label{thm:preimage-of-basic-sets}
Let $f: M^n\to M^n$ be a Smale-Vietoris A-diffeomorphism of closed $n$-manifold $M^n$ and $\mathbb{T}^k\times N=\mathfrak{B}\subset M^n$ the support of Smale skew-mapping $f|_{\mathfrak{B}}=F$, see (\ref{eq:degree-d}). Let $\Omega$ be a basic set of $g: \mathbb{T}^k\to\mathbb{T}^k$ and $\mathfrak{S}=\cap_{l\ge 0}F^l(\mathbb{T}^k\times N)$. Then $\mathfrak{S}\cap p_1^{-1}(\Omega)$ contains a unique basic set $\Omega_{\mathfrak{S}}$ of $f$. Here, $p_1: \mathbb{T}^k\times N\to \mathbb{T}^k$ is the natural projection on the first factor. Moreover,
\begin{enumerate}
  \item If $\Omega$ is a trivial basic set (isolated periodic orbit) of $g$ then $\Omega_{\mathfrak{S}}$ is also trivial basic set.
  \item If $\Omega$ is a nontrivial basic set of $g$ then $\Omega_{\mathfrak{S}}$ is also nontrivial basic set.
  \item If $\Omega$ is a backward $g$-invariant basic set of $g$, $\Omega=g^{-1}(\Omega)$, (hence, $\Omega$ is nontrivial) then
        $\Omega_{\mathfrak{S}}=\mathfrak{S}\cap p_1^{-1}(\Omega)$.
\end{enumerate}
\end{thm}

For $k=1$, when $\mathbb{T}^1=S^1$ is a circle, the following result says that $NW(F)$ contains a unique nontrivial basic set that is either Smale (one-dimensional) solenoid or a nontrivial zero-dimensional basic set.

\begin{thm}\label{thm:descrip-non-wand-set}
Let $f: M^n\to M^n$ be a Smale-Vietoris A-diffeomorphism of closed $n$-manifold $M^n$ and $\mathbb{T}^1\times N=\mathfrak{B}\subset M^n$ the support of Smale
skew-mapping $f|_{\mathfrak{B}}=F$. Then the non-wandering set $NW(F)$ of $F$ belongs to $\mathfrak{S}=\cap_{l\ge 0}F^l(\mathbb{T}^1\times N)$, and $NW(F)$ contains
a unique nontrivial basic set $\Lambda (f)$ that is either
 \begin{itemize}
     \item a one-dimensional expanding attractor, and $\Lambda(f) = \mathfrak{S}$, or
     \item a zero-dimensional basic set, and $NW(F)$ consists of $\Lambda(f)$ and finitely many (nonzero) isolated
           attracting periodic points plus finitely many (possibly, zero) saddle type isolated periodic points of codimension one stable Morse index.
     \end{itemize}
The both possibilities hold.
\end{thm}

It is natural to consider bifurcations from one type of dynamics to another which can be thought of as a destruction (or, a birth) of Smale solenoid. For simplicity, we represent two such bifurcations for $n=3$ and $M^3=S^3$ a 3-sphere. Recall that a diffeomorphism $f: M\to M$ is $\Omega$-\textit{stable} if there is a neighborhood $U(f)$ of $f$ in the space of $C^1$ diffeomorphisms $Diff^1(M)$ such that $f|_{NW(f)}$ conjugate to every $g|_{NW(g)}$ provided $g\in U(f)$.

\begin{thm}\label{thm:Smale-surgery-for-Smale-diffeo}
There is the family of $\Omega$-stable Smale-Vietoris diffeomorphisms $f_{\mu}: S^3\to S^3$, $0\leq\mu\leq 1$, continuously depending on the parameter $\mu$ such that the non-wandering set $NW(f_{\mu})$ of $f_{\mu}$ is the following:
\begin{itemize}
 \item $NW(f_{0})$ consists of a one-dimensional expanding attractor (Smale solenoid attractor) and one-dimensional
       contracting repeller (Smale solenoid repeller);
 \item for $\mu>0$, $NW(f_{\mu})$ consists of two nontrivial zero-dimensional basic sets and finitely many isolated periodic orbits.
\end{itemize}
\end{thm}

\textit{Acknowledgments}. The authors are grateful to V.~Grines, O.~Pochinka and S.~Gonchenko for useful discussions.
Research partially supported by Russian Fund of Fundamental Inves\-ti\-ga\-tions,  13-01-12452 офи-м, 15-01-03687, 13-01-00589. This work was supported by the Basic Research Programs at the National Research University Higher School of Economics ( “Topological methods in the dynamics”,  project 98) in 2016.


\section{Definitions}

A mapping $F: M\times N\to M\times N$ of the type $F(x;y)=\left(g(x);h(x,y)\right)$ is called a \textit{skew-mapping}. One says also a \textit{skew product transformation} over $g$ or simply, a \textit{skew product}. Denote by $End~(M)$ the space of $C^1$ endomorphisms $M\to M$ i.e., the $C^1$ maps of $M$ onto itself. An endomorphism $g$ is \textit{nonsingular} if the Jacobian $|Dg|\neq 0$. This means that $g$ is a local diffeomorphism. In particular, $g$ is a $d$-cover.
In this paper, we consider nonsingular $g\in End~(M)$, $Dg\neq 0$, that are not a diffeomorphism.

Fix $g\in End~(M)$. A point $x\in M$ is said to be \textit{non-wandering} if given any neighborhood $U(x)=U$ of $x$, there is $m\in\mathbb{N}$ such that
$g^m(U)\cap U\neq\emptyset$. Denote by $NW(g)$ the set of non-wandering points. Clearly, $NW(g)$ is a closed set and $g(NW(g))\subset NW(g)$ i.e., $NW(g)$ is a forward $g$-invariant set. The set $\{x_i\}_{-\infty}^{\infty}$ denoted by $O(x_0)$ is called a $g$-\textit{orbit} of $x_0$ if $g(x_i)=x_{i+1}$ for every integer $i$. A subset $\{x_j,x_{j+1},\ldots,x_{j+r}\}\subset O(x_0)$ consisting of a finitely many points of $O(x_0)$ is called a \textit{compact part} of $O(x_0)$.
A $g$-orbit $\{x_i\}_{-\infty}^{\infty}$ is \textit{periodic} if there is an integer $p\geq 0$ such that $g^p(x_i)=x_{i+p}$ for each $i\in\mathbb{Z}$. Certainly, $NW(g)$ contains all periodic $g$-orbits.

The orbit $O(x_0)$ is said to be \textit{hyperbolic} if there is a continuous splitting of the tangent bundle
 $$ \mathbb{T}_{O(x_0)}M=\bigcup_{i=-\infty}^{\infty}\mathbb{T}_{x_i}M=\mathbb{E}^s\bigoplus\mathbb{E}^u=
    \bigcup_{i=-\infty}^{\infty}\mathbb{E}^s_{x_i}\bigoplus\mathbb{E}^u_{x_i} $$
which is preserved by the derivative $Dg$ such that
 $$ ||Dg^m(v)||\leq c\mu^m||v||,\,\, ||Dg^m(w)||\geq c^{-1}\mu^{-m}||w||\quad\mbox{ for }\,\, v\in\mathbb{E}^s,\,\, w\in\mathbb{E}^u,\,\,\forall m\in\mathbb{N} $$
for some constants $c>0$, $0<\mu<1$ and a Riemannian metric on $\mathbb{T}M$. Note that $\mathbb{E}^u({x_0})$ depends on the  negative semi-orbit $\{x_i\}_{i=-\infty}^0$. It may happen that $\mathbb{E}^u({x_0})\neq\mathbb{E}^u({y_0})$ though $x_0=y_0$ but $O(x_0)\neq O(y_0)$. Such a  phenomenon is impossible for $\mathbb{E}^s({x_0})$, it depends only on $x_0$ \cite{Przytycki1977}.

We say that a nonsingular $g\in End~(M)$ \textit{satisfies axiom A}, in short, $f$ is an A-endomorphism if
\begin{itemize}
  \item the periodic $g$-orbits are dense in $NW(g)$ (it follows that $g(NW(g))=NW(g)$);
  \item all $g$-orbits of $NW(g)$ are hyperbolic, and the corresponding splitting of the tangent bundle $\mathbb{T}_{NW(g)}$ depends continuously
        on the compact parts of the $g$-orbits.
\end{itemize}
Recall that Smale's Spectral Decomposition Theorem says that for Axiom A diffeomorphisms the non-wandering set partitions into nonempty closed invariant sets each of which is transitive. Similar theorem for A-endomorphisms was probed in \cite{BlockL75} (Theorem C), \cite{Przytycki1977} (Theorem 3.11 and Proposition 3.13). Thus, if $g$ is a nonsingular A-endomorphism then the non-wandering set $NW(g)$ is the disjoint union $\Omega_1\cup\ldots\cup\Omega_k$ such that each $\Omega_i$ is closed and invariant, $g(\Omega_i)=\Omega_i$, and $\Omega_i$ contains a point whose $g$-orbit is dense in $\Omega_i$. The $\Omega_i$ are called \textit{basic sets}.

Following Williams \cite{Williams67,Williams74}, we introduce an inverse limit for $g: T\to T$ as follows. Put by definition, $
  \prod_g=\{\,(t_0,t_1,\ldots,t_i,\ldots)\in T^{\mathbb{N}}\,\, :\,\, g(t_{i+1})=t_i,\,\, i\geq 0\,\} $
This set is endowed by the product topology of countable factors. This topology has a basis generating by $(\varepsilon,r)$-neighborhoods
\begin{equation}\label{eq:neighborhoods_product}
U=\{\,\, \{x_i\}_0^{\infty}\in\prod_g\,\, : \,\, x_i\in U_{\varepsilon}(t_i),\,\, 0\leq i\leq r\,\, \mbox{ for some }\,\,\varepsilon > 0,\, r\in\mathbb{N}\,\},
 \end{equation}
where $\{t_0,t_1,\ldots,t_i,\ldots\}\in\prod_g$.
Define the shift map
 $$ \hat{g}: \prod_g\to\prod_g, \,\,\,\, \hat{g}(t_0,t_1,\ldots,t_i,\ldots)=\left(g(t_0),t_0,t_1,\ldots,t_i,\ldots\right), \,\,\, (t_0,t_1,\ldots,t_i,\ldots)\in\prod_g. $$
This map $\hat{g}: \prod_g\to\prod_g$ is called the \textit{inverse limit} of $g$ is a homeomorphism \cite{Robinson-book99,Williams74}.

\section{Proofs of main results}\label{s:proof-main-results}\nopagebreak

We denote by $p_1: \mathbb{T}^k\times N\to\mathbb{T}^k$, $p_2: \mathbb{T}^k\times N\to N$ the natural projections $p_1(t,z)=t$ and $p_2(t,z)=z$. A fiber
$\{t\}\times N\stackrel{\rm def}{=}N_t$ of the trivial fiber bundle $p_1$ is called a $t$-leaf. It follows from (\ref{eq:degree-d}) that $F=f|_{\mathfrak{B}}$
takes a $t$-leaf into $g(t)$-leaf.

Let $t\in\mathbb{T}^k$ and $\varepsilon>0$. We denote by $U_{\varepsilon}(t)$ the $\varepsilon$-neighborhood of the point $t$ i.e.,
$U_{\varepsilon}(t)=\{x\in\mathbb{T}^k\, :\, \varrho(x,t)<\varepsilon\}$ where $\varrho$ is a metric on $\mathbb{T}^k$.

The following technical lemma describes the symbolic model of the restriction $f|_{\mathfrak{S}}$. This lemma is a generalization of the similar classical result by Williams \cite{Williams67, Williams74}.

\begin{lm}\label{lm:M-to-sequence} 
Let $f: M^n\to M^n$ be a Smale-Vietoris diffeomorphism of closed $n$-manifold $M^n$ and $\mathbb{T}^k\times N=\mathfrak{B}\subset M^n$ the support of Smale skew-mapping $f|_{\mathfrak{B}}=F$.  Then the restriction $f|_{\mathfrak{S}}$ is conjugate to the inverse limit of the mapping $g: \mathbb{T}^k\to\mathbb{T}^k$, where
$\mathfrak{S}=\cap_{l\ge 0}F^l(\mathbb{T}^k\times N)$.
\end{lm}

\textsl{Proof}. Recall that given any point $t_0\in\mathbb{T}^k$, $g^{-1}(t_0)$ consists of $d$ points, one says $t^1_0$, $t_0^2$, $\ldots$, $t_0^d\in\mathbb{T}^k$.
Since $F$ is a diffeomorphism on its image, the sets $F(N_{t_0^1})$, $\ldots$, $F(N_{t_0^d})$ are pairwise disjoint,
\begin{equation}\label{eq:different-disk}
 F(N_{t_0^i})\cap F(N_{t_0^j}) = \emptyset ,\quad i\neq j,\quad 1\le i,j\le d,
\end{equation}

Now, for the sake of simplicity, we divide the proof into several steps. The end of the proof of each step will be denoted by $\diamondsuit$.

\begin{step}\label{step:sequence-points}
Given any point $p\in\mathfrak{S}$, there is a unique sequence of points $\{t_i\}_{i=0}^{\infty}$, $t_i\in\mathbb{T}^k$, and the corresponding sequence of the leaves $\{N_{t_i}\}_{i=0}^{\infty}$ such that
\begin{itemize}
    \item $p\in\cdots\subset F^i(N_{t_i})\subset F^{i-1}(N_{t_{i-1}})\cdots\subset F(N_{t_1})\subset N_{t_0}$, $p=\cap_{i\ge 0}F^i(N_{t_i})$;
    \item $t_i=g(t_{i+1})$, $i\ge 0$.
\end{itemize}
\end{step}
\textsl{Proof of Step \ref{step:sequence-points}}. Put $t_0=p_1(p)\in\mathbb{T}^k$. Let $g^{-1}(t_0)=\{t^1_0, t_0^2, \ldots, t_0^d\}$. By (\ref{eq:different-disk}), there is a unique $t_0^j$ such that $p\in F(N_{t_0^j})$. Put by definition $t_0^j=t_1$. Note that $F(N_{t_1})\subset N_{t_0}$. Similarly, $g^{-1}(t_1)$ consists of $d$ points $t_1^1$, $t_1^2$, $\ldots$, $t_1^d$. By (\ref{eq:different-disk}), the sets $F(N_{t_1^1})$, $\ldots$, $F(N_{t_1^d})$ are pairwise disjoint. Since $p\in F^2(\mathbb{T}^k\times N)$, there is a unique $t_1^i$ such that
$p\in F^2(N_{t_1^i})$. Put by definition $t_1^i=t_2$. Note that $p\in F^2(N_{t_2})\subset F(N_{t_1})\subset N_{t_0}$. Continuing by this way, one gets the sequences $\{t_i\}_{i=0}^{\infty}$, $\{N_{t_i}\}_{i=0}^{\infty}$ desired. It follows from (\ref{eq:contraction}) that $diam~F^i(N_{t_i})=diam~(F^i(\{t_i\}\times N))\to 0$ as $i\to\infty$. Hence, $p=\cap_{i\ge 0}F^i(N_{t_i})$.
$\diamondsuit$

Let $\hat{g}: \prod_g\to\prod_g$ be the inverse limit of $g: \mathbb{T}^k\to\mathbb{T}^k$ where
$\prod_g=\{\,(t_0,t_1,\ldots,t_i,\ldots)\in\mathbb{T}^{\mathbb{N}}\,\, :\,\, g(t_{i+1})=t_i,\,\, i\geq 0\,\}$.
For a point $p\in\mathfrak{S}$, denote by $P(t_0,t_1,\ldots,t_i,\ldots)$, $t_i\in\mathbb{T}^k$, the sequence due to Step \ref{step:sequence-points}.
Define the mapping
 $$ \theta:\,\, \mathfrak{S}\to\prod_g,\quad p\longmapsto P(t_0,t_1,\ldots,t_i,\ldots),\,\,\, p\in\mathfrak{S}. $$
\begin{step}\label{step:theta-homeo}
The mapping $\theta$ is a homeomorphism.
\end{step}
\textsl{Proof of Step \ref{step:theta-homeo}}. It follows from (\ref{eq:contraction}) that $\theta$ is injective. Since the intersection of nested sequence of closed subsets is non empty, $\theta$ is surjective. One remains to prove that $\theta$ and $\theta^{-1}$ are continuous. Take a neighborhood $U$ of $\theta(p)$, $p\in\mathfrak{S}$. We can assume that $U$ is an $(\varepsilon,r)$-neighborhood (\ref{eq:neighborhoods_product}),
where $\theta(p)=\{t_0,t_1,\ldots,t_i,\ldots\}\in\prod_g$. Moreover, one can assume that $g^{-1}\left(U_{\varepsilon}(t_i)\right)$ consists of $d$ pairwise disjoint domains for every
$0\leq i\leq r$. Recall that $t_i=g(t_{i+1})$, $i\ge 0$. Therefore, $t_{r-j}=g^j(t_r)$ for all $1\le j\le r$. Similarly, $x_{r-j}=g^j(x_r)$,
$1\le j\le r$. Since $g$ is continuous, there exists $0<\delta\le\varepsilon$ such that the inclusion $x_r\in U_{\delta}(t_r)$ implies $x_i\in U_{\varepsilon}(t_i)$ for all $i=0$, $\ldots$, $r$. The restriction $F|_{\mathfrak{S}}: \mathfrak{S}\to\mathfrak{S}$ is a diffeomorphism. Therefore, there is a (relative) neighborhood $U(p)$
of $p$ in $\mathfrak{S}$ such that $ p_1\left(F^{-i}(U(p))\right)\subset U_{\delta}(t_i)\quad\mbox{ for all }\quad 0\leq i\leq r. $
Taking in mind that $g^{-1}\left(U_{\varepsilon}(t_i)\right)$ consists of $d$ pairwise disjoint domains, $0\leq i\leq r$, we see that that
$\theta\left(U(p)\right)\subset U$. Thus, $\theta$ is continuous. Since $\prod_g$ is compact, $\theta^{-1}$ is also continuous.
$\diamondsuit$

\begin{step}\label{step:diagram-commuts}
One holds $\theta\circ F|_{\mathfrak{S}} = \hat{g}\circ\theta|_{\mathfrak{S}}$.
\end{step}
\textsl{Proof of Step \ref{step:diagram-commuts}}. Take $p\in \mathfrak{S}$ and $\theta(p)=\{t_0, t_1,\ldots,t_i,\ldots\}$ where $t_i=g(t_{i+1})$, $i\ge 0$.
By definition of $\hat{g}: \prod_g\to\prod_g$, one holds $ \hat{g}\circ\theta(p) = \hat{g}\left(\theta(p)\right) = \hat{g}\left(\{t_0,t_1,\ldots,t_i,\ldots\}\right)=\{g(t_0),t_0,t_1,\ldots,t_i\ldots\}. $
It follows from (\ref{eq:pavnom}) that $F(p)\in F(\{t_0\}\times N)\subset N_{g(t_0)}$. Hence, by Step \ref{step:sequence-points}, the sequence of points $\{g(t_0),t_0,t_1,\ldots,t_i\ldots\}$ corresponds to $\theta(F(p))$, since
  $$ F(p)=F\left(\cap_{i\ge 0}F^i(\{t_i\}\times N)\right)=\cap_{i\ge 0}F^{i+1}(\{t_i\}\times N)=\cap_{i\ge 0}F^{i+1}(\{t_i\}\times N)\cap N_{g(t_0)}= $$
  $$ = N_{g(t_0)}\cap F(N_{t_0})\cap F^2(N_{t_1})\cap\ldots\cap F^{i+1}(N_{t_i})\cap\ldots. \diamondsuit $$

It follows from Steps \ref{step:theta-homeo}, \ref{step:diagram-commuts} that the mapping $\theta$ is a conjugacy between $F|_{\mathfrak{S}}$ and $\hat{g}$.
Lemma \ref{lm:M-to-sequence} is proved.
$\Box$

To prove Theorem \ref{thm:preimage-of-basic-sets} we need some previous results.

\begin{lm}\label{lm:non-wandering-points}
Let $\overline{t}=\{t_0,t_1,\ldots,t_i,\ldots,\}\in\prod_g$, $g(t_{i+1})=t_i$, $i\geq 0$. Suppose that $t_i\in NW(g)$ for all $i\geq 0$. Then
$\overline{t}\in NW(\hat{g})$ and $\theta^{-1}(\overline{t})\in NW(F)$.
\end{lm}
\textsl{Proof}. Since $\overline{t}=\{t_0,t_1,\ldots,t_i,\ldots,\}=\{g^r(t_r),g^{r-1}(t_r),\ldots,t_r,\ldots\}$, we can take the $(\varepsilon,r)$-neighborhood $V$ \ref{eq:neighborhoods_product} as follows
  $$ V = \{ \, \{g^r(x_r),g^{r-1}(x_r),\ldots,x_r,\ldots\}\,\,\, :\,\,\, g^i(x_r)\in U_{\varepsilon}\left(g^i(t_r)\right),\,0\leq i\leq r\,\}. $$
Since $g$, $g^2$, $\ldots$, $g^r$ are uniformly continuous, there is $0<\delta\leq\varepsilon$ such that $x\in U_{\delta}(y)$ implies $g^i(x)\in U_{\varepsilon}(g^i(y))$ for all $0\leq i\leq r$. By condition, $t_r\in NW(g)$. Hence, there exists $n_0\in\mathbb{N}$ such that
$g^{n_0}\left(V_{\delta}(t_r)\right)\cap V_{\delta}(t_r)\neq\emptyset$.
It follows that there is a point $x_0\in V_{\delta}(t_r)$ such that $g^{n_0}(x_0)\in V_{\delta}(t_r)$.

Take $\overline{x}_0=\{g^r(x_0),g^{r-1}(x_0),\ldots,x_0,\ldots\}\in\prod_g$. Since $x_0\in V_{\delta}(t_r)$, $g^i(x_0)\in U_{\varepsilon}\left(g^i(t_r)\right)$ for
all $0\leq i\leq r$. Therefore, $\overline{x}_0\in V$. Since $g^{n_0}(x_0)\in V_{\delta}(t_r)$, $g^{n_0+i}(x_0)\in U_{\varepsilon}\left(g^i(t_r)\right)$ for all
$0\leq i\leq r$. Therefore,
 $$ \hat{g}^{n_0}(\overline{x}_0)=\{\,g^{n_0+r}(x_0),g^{n_0+r-1}(x_0),\ldots,g^{n_0}(x_0),\ldots\}\in V. $$
As a consequence, $\hat{g}^{n_0}(V)\cap V\neq\emptyset$ and $\overline{t}\in NW(\overline{g})$. A conjugacy map takes a non-wandering set onto non-wandering set.
By Lemma \ref{lm:M-to-sequence}, $\theta^{-1}(\overline{t})\in NW(F)$.
$\Box$

\begin{cor}\label{cor:from-lm:non-wandering-points}
The following qualities hold $ p_1\left[NW(f_{\mathfrak{B}})\right]=p_1\left[NW(F)\right]=NW(g). $
\end{cor}
\textsl{Proof}. Since the projection $p_1$ is continuous, $p_1\left[NW(F)\right]\subset NW(g)$. Take a point $t_0\in NW(g)$. Since $g$ is an A-endomorphism, $g\left[NW(g)\right]=NW(g)$ \cite{BlockL75,Przytycki1977}. Therefore, there is a sequence $t_i\in NW(g)$ such that $g(t_{i+1})=t_i$ for every $i\geq 0$. It follows from Lemma \ref{lm:non-wandering-points} that $\overline{t}=\{t_0,t_1,\ldots,t_i,\ldots,\}\in NW(\hat{g})$ and $\theta^{-1}(\overline{t})\in NW(F)$. By definition of the mapping $\theta$, $\theta^{-1}(\overline{t})\in p^{-1}_1(t_0)$. Hence, $NW(g)\subset p_1\left[NW(F)\right]$.
$\Box$

\begin{lm}\label{lm:sequence-consists-of-nonwandering}
Let $(t_0,z_0)\in\mathfrak{S}$ be a non-wandering point of $f$, and $\theta(t_0,z_0)=\{t_i\}_{i\geq 0}$. Then $t_i\in NW(g)$ for all $i\geq 0$.
\end{lm}
\textsl{Proof}. According to Corollary \ref{cor:from-lm:non-wandering-points}, $p_1\left[NW(f_{\mathfrak{B}})\right]=p_1\left[NW(F)\right]=NW(g)$. Therefore,
$t_0\in NW(g)$. Since $F_{\mathfrak{S}}: \mathfrak{S}\to\mathfrak{S}$ is a diffeomorphism, $F^{-1}\left(NW(F)\right)=NW(F)$ and
$F^{-1}(t_0,z_0)=(t_1,z_1)\in NW(F)=NW(f_{\mathfrak{B}})$. Hence, $t_1\in NW(g)$ by Step 1. Continuing this way, one gets that $t_i\in NW(g)$ for all $i\geq 0$.
$\Box$

\begin{cor}\label{cor1:from-lm-sequence-consists-of-nonwandering}
Let $(t_0,z_0)\in\mathfrak{S}$ be a non-wandering point of $f$, and $\theta(t_0,z_0)=\{t_i\}_{i\geq 0}$. Suppose that $t_0$ belongs to a basic set $\Omega$ of $g$. Then $t_i\in\Omega$ for all $i\geq 0$.
\end{cor}
\textsl{Proof}. By Lemma \ref{lm:sequence-consists-of-nonwandering}, $t_i\in NW(g)$ for all $i\geq 0$. Since $\Omega$ is forward $g$-invariant, $t_i\in\Omega$ for all $i\geq 0$.
$\Box$

\begin{lm}\label{lm:same-basic-set}
Let $\Omega$ be a nontrivial basic set of $g$, and $t_0\in\Omega$. Suppose that two points $(t_0,z_1)$, $(t_0,z_2)\in\mathfrak{S}$ are non-wandering under $f$. Then the both $(t_0,z_1)$ and $(t_0,z_2)$ belong to the same basic set of $f$.
\end{lm}
\textsl{Proof}. Denote by $\Omega_j$ the basic set of $F$ containing the point $(t_0,z_j)$, $j=1,2$. Clearly, $\Omega_j\subset\mathfrak{S}$. We have to prove that $\Omega_1=\Omega_2$. It is sufficient to show that there is a non-wandering point $q\in NW(F)$ such that each point $(t_0,z_1)$ and $(t_0,z_2)$ belongs to the $\omega$-limit set of $q$.

Let $\overline{t}_j=\theta(t_0,z_j)=\{t_0,t_1^{(j)},\ldots,t_i^{(j)},\ldots\}$, $j=1,2$. By Corollary \ref{cor1:from-lm-sequence-consists-of-nonwandering}, $t_i^{(j)}\in\Omega$ for all $i\geq 0$, $j=1,2$. Since the basic set $\Omega$ is transitive, there is a point $x_0\in\Omega$ such that its positive semi-orbit $O^+_g(x_0)$ is dense in $\Omega$, $clos~\left(O^+_g(x_0)\right)=\Omega$.

It follows from Corollary \ref{cor:from-lm:non-wandering-points} that there is a point $\overline{x}_0=\{x_0,x_1,\ldots,x_i,\ldots\}\in\prod_g$ such that $x_i\in\Omega$ for all $i\geq 0$. Take arbitrary $(\varepsilon,r)$-neighborhood $U(\overline{t}_1)$ of $\overline{t}_1$. Since $g$, $g^2$, $\ldots$, $g^r$ are uniformly continuous, there exists $\delta>0$ such that the inequality
$x\in U_{\delta}(y)$ implies $g^i(x)\in U_{\varepsilon}(y)$ for all $0\leq i\leq r$. Because of the semi-orbit $O^+_g(x_0)$ is dense in $\Omega$, there is $n_0\in\mathbb{N}$ such that $g^{n_0}(x_0)\in U_{\delta}(t^{(1)})$. Hence, $\hat{g}^{n_0}(\overline{x}_0)\in U(\overline{t}_1)$. Therefore,  $\overline{t}_1=\theta(t_0,z_1)$ belongs to the $\omega$-limit set of $\overline{x}_0$. Similarly, one can prove that  $\overline{t}_2=\theta(t_0,z_2)$ belongs to the $\omega$-limit set of $\overline{x}_0$ as well. Since $\theta$ is a conjugacy mapping, the points $(t_0,z_1)=\theta^{-1}(\overline{t}_1)$ and $(t_0,z_2)=\theta^{-1}(\overline{t}_2)$ belongs to the $\omega$-limit set of the point $q=\theta^{-1}(\overline{x}_0)\in NW(F)$.
$\Box$

\textbf{\textsl{Proof of Theorem \ref{thm:preimage-of-basic-sets}.}} We know that $p_1\left[NW(F)\right]=NW(g)$. Hence, $\mathfrak{S}\cap p_1^{-1}(\Omega)$ contains basic sets
of $f$. Suppose that $\Omega$ is trivial i.e., $\Omega$ is an isolated periodic orbit
 $$ \Omega=Orb_g(q)=\{q,g(q),\ldots,g^{p-1}(q),g^p(q)=q\},\,\,\,\mbox{ where }\,\,\, q\in\mathbb{T}^k\,\, \mbox{ and }\,\, p\in\mathbb{N}\,\,
    \mbox{ is a period of }\,\, q. $$
By definition of Smale skew-mapping, the restriction of $F=f|_{\mathfrak{B}}$ on the second factor $N$ is the uniformly attracting embedding. Therefore,
 $$ N_{q}\supset f^p(N_{q})\supset\cdots\supset f^{mp}(N_{q})\supset\cdots\, \mbox{ and the intersection }\,\, \bigcap_{m\geq 0}f^{mp}(N_{q})\,\,
    \mbox{ is a unique point, say } Q. $$
Similarly, $\cap_{m\geq 0}f^{mp}(N_{g^i(q)})$ is a unique point $f^i(Q)$ for every $0\leq i\leq p-1$. It follows from (\ref{eq:degree-d}), that $\{Q,f(Q),\ldots,f^{p-1}(Q),f^p(Q)=Q\}$ is an isolated periodic orbit $Orb_f(Q)$ such that $NW(F)\cap p_1^{-1}(\Omega)=Orb_f(Q)$. Therefore, $Orb_f(Q)=\Omega_{\mathfrak{S}}$ is a unique basic set of $F$ that belongs to
$\mathfrak{S}\cap p_1^{-1}(\Omega)$.

Let $\Omega$ be a nontrivial basic set. It follows from Lemma \ref{lm:same-basic-set} that all basic set of $F$ that is contained in $\mathfrak{S}\cap p_1^{-1}(\Omega)$ are coincide. Hence, $\Omega_{\mathfrak{S}}$ is a unique nontrivial basic set of $f$ contained in $\mathfrak{S}\cap p_1^{-1}(\Omega)$.

Now let $\Omega$ be a backward $g$-invariant basic set of $g$. Note that the equality $\Omega=g^{-1}(\Omega)$ implies that $\Omega$ cannot be a trivial basic set, since $g$ is a $d$-cover, $d\geq 2$. It follows from Lemma \ref{lm:non-wandering-points} that every point of $\mathfrak{S}\cap p_1^{-1}(\Omega)$ is a non-wandering point of $f$. By Lemma \ref{lm:same-basic-set}, $\mathfrak{S}\cap p_1^{-1}(\Omega)$ is a unique basic set. Theorem \ref{thm:preimage-of-basic-sets} is proved.
$\Box$

\textsl{Example.} Let us consider three endomorphism $g_i: \mathbb{T}^2\to\mathbb{T}^2$, $i=1,2,3$, that are 2-covers. $g_1$ is defined by the matrix
$\left( \begin{array}{cc}
  3 & 1 \\
  1 & 1
\end{array} \right).$  Clearly, $g_1$ is an expanding A-endomorphism, and $\mathbb{T}^2$ is a unique basic set of $g_1$. The corresponding diffeomorphism $f$ has a unique basic set, say $\Omega_1$, thai is locally homeomorphic to the product of $\mathbb{R}^2$ and Cantor set. Thus, $\Omega_1$ is 2-dimensional.

Now, let us consider the case when $\mathbb{T}^1=S^1$ is a circle, and $d$-cover $g: \mathbb{T}^1\to\mathbb{T}^1$ is a nonsingular endomorphism of $S^1$. The crucial step of the proof of Theorem  \ref{thm:descrip-non-wand-set} is the following result.

\begin{lm}\label{lm:non-wandering-set-a-endomorphism}
Let $g: \mathbb{T}^1\to\mathbb{T}^1$ be a nonsingular A-endomorphism, and $NW(g)$ a non-wandering set of $g$. Then $NW(g)$ is either $\mathbb{T}^1$ or $NW(g)$ is
the union of the Cantor type set $\Sigma$ and finitely many (nonzero) isolated attracting periodic orbits plus finitely many (possibly, zero) repelling isolated
periodic orbits. Moreover, in the last case, $\Sigma$ is backward $g$-invariant.
\end{lm}
\textsl{Proof}. Suppose that $NW(g)\neq\mathbb{T}^1$. By \cite{Shub69}, $g$ is semi-conjugate to the expanding linear mapping $E_d$, $E_d(t)=dt\,\, mod\, 1$, i.e.,
there is a continuous map $h: \mathbb{T}^1\to\mathbb{T}^1$ such that $g\circ h=h\circ E_d$. Moreover, $h$ is monotone \cite{Melo-Strien-book-1993}. As a consequence,
given any point $t\in\mathbb{T}^1$, $h^{-1}(t)$ is either a point or a closed segment. Since $NW(g)\neq\mathbb{T}^1$, $h$ is not a homeomorphism. Hence, there are
points $t\in\mathbb{T}^1$ for which $h^{-1}(t)$ is a (nontrivial) closed segment. Denote the set of such points by $\chi$. The set $\chi$ is countable and invariant
under $E_d$, $E_d(\chi)=E_d^{-1}(\chi)=\chi$ \cite{ABZ,Melo-Strien-book-1993}. Therefore, $h^{-1}(\chi)$ is also invariant under $g$. As a consequence, $\Sigma=\mathbb{T}^1\setminus clos~(h^{-1}(\chi))$ is the Cantor set consisting on non-wandering points of $g$. Moreover, $\Sigma$ is invariant under $g$
(in particular, backward $g$-invariant). It follows from \cite{Nitecki70} that the part of $NW(g)$ that different from $\Sigma$ consists of finitely many (nonzero) isolated attracting periodic orbits and finitely many (possibly, zero) repelling isolated periodic orbits.
$\Box$

Now, Theorem \ref{thm:descrip-non-wand-set} except the realization part immediately follows from Theorem \ref{thm:preimage-of-basic-sets} and Lemma \ref{lm:non-wandering-set-a-endomorphism}. It remains to construct a Smale-Vietoris A-diffeomorphism the non-wandering set of whose consists of a nontrivial zero-dimensional basic set and a finitely many (nonzero) isolated periodic orbits. It follows from \cite{Bothe83,JiagNiWang2004} for $n=3$ and \cite{BlockL75,Williams67} for $n\geq 4$ that it is sufficient to construct Smale skew-mapping $F: S^1\times D^{n-1}\to S^1\times D^{n-1}$ with the non-wandering set desired because of Smale skew-mapping can be extended to a diffeomorphism of some closed $n$-mannifold. Moreover, according to Robinson-Williams \cite{RobinsonWilliams76} construction of classical Smale solenoid, we can suppose $n=3$.

Let $g: S^1\to S^1$ be a $C^{\infty}$ nonsingular A-endomorphism that is a $d$-cover ($d\ge 2$) with the non-wandering set $NW(g)$ consisting of a unique attracting fixed point $x_0$ and a Cantor set $\Omega$. Moreover, one can assume that $Dg|_{\Omega}=2d-1$, $Dg(x_0)=\lambda < 1$ where $\lambda$ will be specified below. Such endomorphism was constructed by Shub \cite{Shub69}. Hirsch \cite{Hirsch74} has noticed that such endomorphism can be smoothed to be analytical. Now, the circle $S^1$ is endowed with the parameter inducing by the natural projection $[0;1]\to [0;1]/(0\sim 1)=S^1$. We can assume that the restriction $g|{[0;\frac{1}{2}]}$ is a diffeomorphism $[0;\frac{1}{2}]\to [0;\frac{1}{2}]$ with the attracting fixed point $x_0=\frac{1}{4}$ and two repelling fixed points $0$, $\frac{1}{2}$. Without loss of generality, one can also assume that $g|_{[\frac{1}{2};1]}(x)=(2d-1)x\mod 1$. By construction, $\cup_{n\ge 0}g^{-n}_d\left(0;\frac{1}{2}\right)$ is the stable manifold $W^s(x_0)$ of $x_0$, and $\Omega = S^1\setminus W^s(x_0)$ is Cantor set belonging to $NW(g)$. Clearly, given any $y\in S^1$, $\min_{t_k,t_j}\{|t_k-t_j|=\frac{1}{2d-1}$ where
$t_k\neq t_j$ and $g(t_k)=g(t_j)=y$. We take $0 < \lambda < \frac{1}{4}\sin \frac{\pi}{2d-1}$.
After this specification, we denote $g$ by $g_d$. Put by definition
\begin{equation}\label{eq:particular-map-F}
    F(t,z)=\left(g_d(t),\, \lambda z + \frac{1}{2}\exp{2\pi it}\right),\,\,\, F: \mathcal{B}=S^1\times D^2\to\mathcal{B},
\end{equation}
where $D^2\subset\mathbb{R}^2$ is the unit disk, and $z=x+iy$, and $\mathcal{B}$ is a support of Smale skew-mapping. Since $\lambda < \frac{1}{4}$,
$F(\mathcal{B})\subset int~\mathcal{B}$. The Jacobian of $F$ equals
\begin{equation}\label{eq:jacobian}
    DF(t,z)=
\left(
  \begin{array}{cc}
    Dg_d(t) & 0 \\
    \pi i\exp{2\pi it} & \lambda Id_2 \\
  \end{array}
\right),
\end{equation}
where $Id_2$ is the identity matrix on $\mathbb{C}$ or $\mathbb{R}^2$. Since $Dg_d>0$ and $\lambda > 0$, $F$ is a local diffeomorphism. It follows from
$\lambda < \frac{1}{4}\sin \frac{\pi}{2d-1}$ that $F$ is a (global) diffeomorphism on its image.

Since $g_d$ is an A-endomorphism, the periodic points of $g_d$ are dense in $NW(g_d)$. By Lemma \ref{lm:non-wandering-points}, the periodic points of $F$ are dense in $NW(F)$. Thus, it remains to prove the $NW(F)$ has a hyperbolic structure. We follow \cite{Robinson-book99}, Proposition 8.7.5. Clearly, the tangent bundle $T(\mathcal{B})=T(S^1\times D^2)$ is the sum $T(\mathcal{B})=T(S^1)\oplus T(D^2)$, and the fiber $T_{(t,z)}(\mathcal{B})$ at each point $(t,z)\in\mathcal{B}$ is the sum of one-dimensional and two-dimensional tangent spaces $T_t(S^1)=\mathbb{E}^1\cong\mathbb{R}$, $T_z(T^2)=\mathbb{E}^2\cong\mathbb{R}^2$ respectively. It follows from (\ref{eq:jacobian}) that $\mathbb{E}^2$ is invariant under $DF$:
 $$ DF_p{\vec 0\choose \vec v_{23}} = {\vec 0\choose \lambda\vec v_{23}},\quad \vec v_{23}\in \mathbb{E}^2. $$
Moreover, since $|\lambda|<1$, $\mathbb{E}^2$ is the stable bundle, $E^s=\mathbb{E}^2$.

Take $q=(t,z)\in NW(F)$. Then $p_1(q)=t\in NW(g_d)$. If $t=x_0$, then $q$ is a hyperbolic (attractive) fixed point of $F$. For $t\in\Omega$, we consider the cones
 $$ C^u_q = \left\{{\vec v_1\choose \vec v_{23}}\, :\, \vec v_1\in T_t(S^1),\, \vec v_{23}\in \mathbb{E}^2_z,
    \quad |\vec v_1|\ge \frac{2d-1}{4}|\vec v_{23}|\right\}\subset T(\mathcal{B})=\mathbb{E}^1\oplus \mathbb{E}^2. $$
For ${\vec v_1\choose \vec v_{23}}\in C^u$, it follows from (\ref{eq:jacobian}) that
 $$ DF{\vec v_1\choose \vec v_{23}} = {\vec v_1^{\prime}\choose \vec v_{23}^{\prime}} =
    {(2d-1)\vec v_1\choose \pi i\vec v_1\exp{2\pi it}+\lambda\vec v_{23}}. $$
Hence, $|\vec v_{23}|\le |\pi i\exp{2\pi it}\vec v_1|+\lambda|\vec v_{23}|=\pi |\vec v_1|+\lambda|\vec v_{23}|$. Taking in mind $\lambda\le\frac{1}{4}$, one gets
 $$ |\vec v_1^{\prime}|=(2d-1)|\vec v_1|=\frac{2d-1}{4}(4|\vec v_1|)\ge
    \frac{2d-1}{4}\left(\pi|\vec v_1|+\frac{1}{2}|\vec v_1|\right)\ge $$
 $$ \ge\frac{2d-1}{4}\left(\pi|\vec v_1|+\frac{2d-1}{8}|\vec v_{23}|\right)\ge
    \frac{2d-1}{4}\left(\pi|\vec v_1|+\lambda|\vec v_{23}|\right)\ge\frac{2d-1}{4}|\vec v_{23}^{\prime}|, $$
since $\frac{2d-1}{8}\ge\frac{1}{4}$. Therefore, ${\vec v_1^{\prime}\choose \vec v_{23}^{\prime}}\in C^u_{F(q)}$ and $DF(C^u_q)\subset C^u_{F(q)}$. As a consequence,
$$ DF^k(C^u_{F^{-k}(q)})\subset DF^{k-1}(C^u_{F^{-k+1}(q)})\subset\cdots\subset DF(C^u_{F^{-1}(q)})\subset C^u_q
    \quad\mbox{for any}\quad k\in\mathbb{N}. $$
To prove that the intersection of this nested cones is a line, take
 $$ {\vec v_1\choose \vec v_{23}},\, {\vec w_1\choose \vec w_{23}}\in C^u_{F^{-k}(q)},\quad
    {\vec v_1^k\choose \vec v_{23}^k}=DF^k{\vec v_1\choose \vec v_{23}},\,
    {\vec w_1^k\choose \vec w_{23}^k}=DF^k{\vec w_1\choose \vec w_{23}}. $$
Put by definition, $|\vec v_1^j|=v_1^j$, $|\vec w_1^j|=w_1^j$, $\vec v_1=(v_1,0)$, $\vec w_1=(w_1,0)$, $v_1>0$, $w_1>0$. Then
 $$ \left|\frac{\vec v_{23}^1}{v_1^1}-\frac{\vec w_{23}^1}{w_1^1}\right| =
    \left|\frac{\pi i\vec v_1\exp{2\pi it}+\lambda\vec v_{23}}{(2d-1)v_1} -
    \frac{\pi i\vec w_1\exp{2\pi it}+\lambda\vec w_{23}1}{w_1}\right| = $$
 $$ \left|\frac{\pi i\exp{2\pi it}(w_1\vec v_1-v_1\vec w_1)}{(2d-1)v_1w_1} +
    \frac{\lambda}{2d-1}\left(\frac{\vec v_{23}}{v_1}-\frac{\vec w_{23}}{w_1}\right)\right| =
    \frac{\lambda}{2d-1}\left|\frac{\vec v_{23}}{v_1}-\frac{\vec w_{23}}{w_1}\right|, $$
since $w_1\vec v_1-v_1\vec w_1=|\vec w_1|\vec v_1-|\vec v_1|\vec w_1=0$. Therefore,
 $$ \left|\frac{\vec v_{23}^k}{v_1^k}-\frac{\vec w_{23}^k}{w_1^k}\right| =
    \left(\frac{\lambda}{2d-1}\right)^k\left|\frac{\vec v_{23}}{v_1}-\frac{\vec w_{23}}{w_1}\right|, $$
which goes to 0 as $k$ goes to $\infty$. Since the difference of slopes goes to 0, the cones converge to a line, say $\mathbb{E}^u$. The calculation gives that the restriction of the derivative $DF$ on $\mathbb{E}^u$ is an expansion.
$\Box$

Taking in mind the realization part of the proof of Theorem \ref{thm:descrip-non-wand-set}, we see that it is sufficient to construct the corresponding family of $d$-endomorphisms $S^1\to S^1$, $d\geq 2$. First, we represent the two parameter family of circle endomorphisms $f_{\varepsilon,\delta}$ continuously depending on the parameters $\varepsilon\in (0;1)$ and $\delta\in [0;\frac{1}{4})$.

Let $U_{\delta}(x)$ be the bump-function such that
\begin{itemize}
  \item $U_{\delta}(x)=1$ for $x\in\left[-\frac{\delta}{2};+\frac{\delta}{2}\right]$, $0<\delta\leq\frac{1}{4}$;
  \item $U_{\delta}(x)=0$ for $|x|\geq\delta$;
  \item $U_{\delta}'(x)\geq 0$ for $x\in\left[-\delta;-\frac{\delta}{2}\right]$,$\quad$ and $\quad$
        $U_{\delta}'(x)\leq 0$ for $x\in\left[\frac{\delta}{2};\delta\right]$.
\end{itemize}

\begin{lm}\label{lm:short-bifurcation-d-endomorphisms}
Let
 $$ f_{\varepsilon,\delta}(x)=\left\{
 \begin{array}{ccc}
   dx+(-d+\varepsilon)xU_{\delta}(x)\,\,\, mod \,\, 1 & \mbox{ for } & \varepsilon\in (0;1),\,\,\delta\in (0;\frac{1}{4}) \\
   dx\,\,\, mod \,\, 1 & \mbox{ for } & \varepsilon=0,\,\, \delta=0
 \end{array}\right.
 $$
Then $f_{\varepsilon,\delta}$ is a structurally stable nonsingular circle $d$-endomorphism such that the non-wandering set $NW(f_{\varepsilon,\delta})$ is the union of
a unique hyperbolic attracting point $x=0$ and a Cantor set provided $\varepsilon\neq 0$ and $\delta\neq 0$. Moreover, $NW(f_{0,0})=S^1$. In addition, $f_{\varepsilon,\delta}\to E_d$ as $\varepsilon\neq 0$ is fixed and $\delta\to 0$ in the $C^0$ topology.
\end{lm}
\textsl{Proof}. For $\varepsilon\neq 0$ and $\delta\neq 0$, we see
 $$ f_{\varepsilon,\delta}'(x)=d+(-d+\varepsilon)\left[xU_{\delta}(x)'+U_{\delta}(x)\right]=
    d+(-d+\varepsilon)xU_{\delta}(x)'+(-d+\varepsilon)U_{\delta}(x). $$
Clearly, $d+(-d+\varepsilon)U_{\delta}(x)\geq\varepsilon$. Since $xU_{\delta}(x)'\leq 0$, $f_{\varepsilon,\delta}'(x)\geq\varepsilon$. Because of outside
of the $\delta$-neighborhood $V_{\delta}(0)$ of $x_0=0$ the mapping $f_{\varepsilon,\delta}$ coincides with the linear $d$-endomorphism $E_d(x)=dx\,\, mod \,\, 1$, $f_{\varepsilon,\delta}$ is a nonsingular $d$-endomorphism. Since $f_{\varepsilon,\delta}'(0)=\varepsilon\in (0;1)$, $x=0$ is a hyperbolic
attracting point. Solving the equation $dx+(-d+\varepsilon)xU_{\delta}(x)=x$, one gets two fixed points $\pm x_*\in V_{\delta}(0)$ such that
$U_{\delta}(\pm x_*)=\frac{d-1}{d-\varepsilon}$, where $\frac{\delta}{2}<x_*<\delta$. Moreover, the $\omega$-limit set of any point from $(-x_*;x_*)$ is $x_0=0$. Hence, $NW(f_{\varepsilon,\delta})$ equals
 $$ NW(f_{\varepsilon,\delta})=\{x_0\}\bigcup\left(S^1\setminus\cup_{k\geq 0}f^{-k}_{\varepsilon,\delta}(-x_*;x_*)\right), $$
where $C=S^1\setminus\cup_{k\geq 0}f^{-k}_{\varepsilon,\delta}(-x_*;x_*)$ is Cantor set. For any $x\in C$, one can prove that
 $$ f_{\varepsilon,\delta}'(x)=d+(-d+\varepsilon)xU_{\delta}'(x)+(-d+\varepsilon)U_{\delta}(x)
    \geq d+(-d+\varepsilon)U_{\delta}(x_*)+(-d+\varepsilon)xU_{\delta}'(x)= $$
 $$ = 1+(-d+\varepsilon)xU_{\delta}'(x)>1. $$
It follows from \cite{Nitecki70} that $f_{\varepsilon,\delta}$ is structurally stable. At last, for $x\in V_{\delta}(0)$, one gets
 $$ \left| f_{\varepsilon,\delta}(x)-E_d(x)\right|=\left|(-d+\varepsilon)xU_{\delta}(x)\right|\leq\delta d\to 0\,\,\,\mbox{ as }\,\,\, \delta\to 0. $$
As a consequence, $f_{\varepsilon,\delta}\to E_d$ as $\delta\to 0$ in the $C^0$ topology.
$\Box$

Taking in mind Lemma \ref{lm:short-bifurcation-d-endomorphisms} and using the technics developed in \cite{JiagNiWang2004} (see also \cite{BlockL75,Bothe83,{JimingMaBinYu2007},Williams67}), one can prove Theorem \ref{thm:Smale-surgery-for-Smale-diffeo}.

\renewcommand{\refname}{References}

\noindent National Research University Higher School of Economics, \, Nizhny Novgorod Academy MVD of Russia

\noindent \textit{E-mail:} zhuzhoma@mail.ru,\, nisaenkova@mail.ru

\end{document}